\documentclass[11pt]{article}
\usepackage{amsfonts}
\usepackage{latexsym,tikz}
\usepackage{amsmath}
 
\title{Multiple Dedekind Zeta Functions}
\author{Ivan Horozov}
\usepackage{color}
\date{}

\newtheorem{theorem}{Theorem} %[section]
\newtheorem{lemma}[theorem]{Lemma}
\newtheorem{proposition}[theorem]{Proposition}
\newtheorem{corollary}[theorem]{Corollary}

\def \N {{\mathbb N}}
\def \Z {{\mathbb Z}}
\def \Q {{\mathbb Q}}

\def \P {{\mathbb P}}

\begin{document}

%%%%%%%%%%%%%%%%%%%%%%%%%%%%%%%%%%%%%%%%%%%%%%%%%%%%%%%%%%%%%%%%%%%%%%%
%%%%%%
\title{Periods of Mixed Tate Motives over Real Quadratic Number Rings}
\author{Ivan Horozov}
\maketitle

\begin{abstract} Recently, the author defined multiple Dedekind zeta values \cite{MDZF} associated to a number $K$ field and a cone $C$.
In this paper we construct explicitly non-trivial examples of mixed Tate motives over the ring of integers in $K$, for a real quadratic number field $K$ and a particular cone C. The period of such a motive is a multiple Dedekind zeta values at $(s_1,s_2)=(1,2)$, associated to the pair $(K;C)$, times a nonzero element of $K$.
%We prove that over any real quadratic field $K$ there is a cone so that the corresponding multiple Dedekind zeta value at $(s_1,s_2)=(1,2)$ is a mixed Tate period over the ring of integers in the real quadratic field $K$.
\end{abstract}
%%%%%%%%%%%%%%%%%%%%%%%%%%%%%%%%%%%



%%%%%%%%%%%%%%%%%%%
%\tableofcontents
%\setcounter{section}{-1}
%%%%%%%%%%%

\section{Introduction}
The Riemann zeta function 
\[\zeta(s)=\sum_{n>0}\frac{1}{n^s}\]
is widely used in number theory, algebraic geometry and quantum field theory.
Euler's multiple zeta values
\[\zeta(s_1,\dots,s_m)=\sum_{0<n_1<\dots<n_m}\frac{1}{n_1^{s_1}\dots n_m^{s_m}},\]
where $s_1,\dots,s_m$ are positive integers and $s_m\geq 2$,
appear  as values of some Feynman amplitudes, and in algebraic geometry, as periods of mixed Tate motives over $Spec(\Z)$ (see \cite{GM}, \cite{DG}, \cite{B1}, \cite{KZ}).

Dedekind zeta values
\[\zeta_K(s)=\sum_{\mathfrak{a}\neq(0)}\frac{1}{N(\mathfrak{a})^s},\]
are a generalization of the Riemann zeta function to a number field $K$.
In some Feynman amplitudes one of the summands is $\log(1+\sqrt{2})$ or $\log\left(\frac{1+\sqrt{5}}{2}\right)$. These values are essentially the residues at $s=1$ of Dedekind zeta functions over $\Q(\sqrt{2})$ and over $\Q(\sqrt{5})$, respectively. 
For $s=2,3,4,\dots$ the values $\zeta_K(s)$ are periods of mixed Tate motives over the ring of algebraic integers in $K$ with ramification only at the discriminant of $K$ (see \cite{B2}).

In \cite{MDZF}, the author has constructed multiple Dedekind zeta values, which are a generalization of Euler's multiple zeta values to number fields in the same way as Dedekind zeta values generalizes Riemann zeta values. For a quadratic number field $K$, the key examples of multiple Dedekind zeta values are
\begin{multline}
\zeta_{K;C}(s_1,\dots,s_1;\dots;s_m,\dots,s_m)=\\
\sum_{\alpha_1,\dots,\alpha_m\in C}
\frac{1}{N(\alpha_1)^{s_1}N(\alpha_1+\alpha_2)^{s_2}\cdots N(\alpha_1+\dots+\alpha_m)^{s_m}},
\end{multline}
where $s_1,\dots,s_m$ are positive integers and $s_m\geq 2$ and
$C$ is a cone generated by a totally positive unit $\beta$ in $K$ and $1$, defined by
\[C=\N\{1,\beta\}=\{\gamma\in K\,|\, \gamma=a+b\beta,\mbox{ for positive integers } a \mbox{ and }b\}.\]
Similar types of cones were considered by Zagier in \cite{Z1} and \cite{Z2}.

In \cite{MDZF}, the author has proven that multiple Dedekind zeta values can be interpolated to multiple Dedekind zeta functions, which have meromorphic continuation to all complex values of the variables $s_1,\dots,s_m$.

In this paper we prove the following theorem.
\begin{theorem}
\label{thm}
Let $K$ be a real quadratic field, and let $C$ be a cone generated by a totally positive unit $\beta$ in $K$ and  $1$.
Then the multiple Dedekind zeta values
\[(\beta_2-\beta_1)^{3}\zeta_{K;C}(1,2)\]
is a period of a mixed Tate motive over the ring of integers in $K$. 
In particular, it is unramified over the primes dividing the discriminant $\sqrt{D}$.
\end{theorem}

{\bf{Remark:}} The proof of the Theorem can easily be generalized to all 
\[(\beta_2-\beta_1)^{s_1+\cdots+s_m}\zeta_{K;C}(s_1,\dots,s_m)\] 
for the same cone $C$. The details for the general case will be completed in a sequel to this paper. The choice of considering $\zeta_{K;C}(1,2)$ in this paper is two-fold. First, this is among the simplest non-trivial example of a multiple Dedekind zeta value. Second, for any other (multple) Dedekind zeta value, the proof of the corresponding statement is essentially the same.
%%%%%
\section{Background}

\subsection{Multiple zeta values}
The Riemann zeta function at the value $s=2$ can be expressed in term of an iterated integral in the following way

\begin{multline*}
\int_0^1\left(\int_0^y\frac{dx}{1-x}\right)\frac{dy}{y}
=\int_0^1\left(\int_0^y(1+x+x^2+x^3\dots)dx\right)\frac{dy}{y}\\
=\int_0^1\left(y+\frac{y^2}{2}+\frac{y^3}{3}+\frac{y^4}{4}+\dots \right)\frac{dy}{y}
=\left. y+\frac{y^2}{2^2}+\frac{y^3}{3^2}+\frac{y^4}{4^2}\cdots\right|_{y=0}^{y=1}\\
=1+\frac{1}{2^2}+\frac{1}{3^2}+\frac{1}{4^2}\dots
=\zeta(2).
\end{multline*}

Let us examine the domain of integration of the iterated integral. Note that $0<x<y$ and $0<y<1$. We can put both inequalities together. Then we obtain the domain $0<x<y<1$, which is a simplex. Thus, we can express the iterated integral as
\[
\zeta(2)=\int_0^1\left(\int_0^y\frac{dx}{1-x}\right)\frac{dy}{y}=\int_{0<x<y<1}\frac{dx}{1-x}\wedge\frac{dy}{y}.
\]

Moreover, Goncharov and Manin \cite{GM} have expressed all multiple zeta values as periods of motives  related to the moduli space of curves of genus zero with $n+3$ marked points, $\mathcal{M}_{0,n+3}$. In particular, $\zeta(2)$ can be expressed as a period of the motive $H^2(\overline{\mathcal{M}}_{0,5}-A,B-A\cap B)$ by pairing of $[\Omega_A]\in Gr^W_4H^2(\overline{\mathcal{M}}_{0,5}-A)$ for $\Omega_A=\frac{dx}{1-x}\wedge\frac{dy}{y}$, 
with $[\Delta_B]\in \left(Gr^W_0H^2(\overline{\mathcal{M}}_{0,5}-B)\right)^\vee$. The Deligne-Mumford compactification $\overline{\mathcal{M}}_{0,5}$ of the moduli space $\mathcal{M}_{0,5}$ can be obtained by three blow-ups of $\P^1\times\P^1$ at the points $(0,0)$, $(1,1)$ and $(\infty,\infty)$. 
Let us name the exceptional divisors at the three points by $E_0$, $E_1$ and $E_\infty$, respectively. Then
$A=(x=1)\cup(y=0)\cup (x=\infty)\cup (y=\infty)\cup E_\infty$ and $B=(x=0)\cup(x=y)\cup(y=1)\cup E_0\cup E_1$.

Similarly, one can express $\zeta(3)$ and $\zeta(1,2)$ as iterated integrals
{\small \begin{gather*}
\zeta(3)=\int_0^1\left(\int_0^z\left(\int_0^y\frac{dx}{1-x}\right)\frac{dy}{y}\right)\frac{dz}{z}=\int_{0<x<y<z<1}\frac{dx}{1-x}\wedge\frac{dy}{y}\wedge\frac{dz}{z},\\
\zeta(1,2)=\int_0^1\left(\int_0^z\left(\int_0^y\frac{dx}{1-x}\right)\frac{dy}{1-y}\right)\frac{dz}{z}=
\int_{0<x<y<z<1}\frac{dx}{1-x}\wedge\frac{dy}{1-y}\wedge\frac{dz}{z}.
\end{gather*}}
Again, $\zeta(3)$ and $\zeta(1,2)$ can be expressed as periods of motives related to $\mathcal{M}_{0,6}$.
In the same paper, Goncharov and Manin prove that the motives associated to multiple zeta values (MZVs) are mixed Tate motives unramified over $Spec(\Z)$.

A few years later, Francis Brown \cite{B1} proved that periods of mixed Tate motives unramified over $Spec(\Z)$ can be expressed as a $\Q$-linear combination of MZVs times an integer power of $2\pi i$.

%\subsection{Dedekind zeta function}
%Dedekind generalized the Riemann zeta function to number fields. Let us recall part of this construction.
%Let $f(x)=x^n+c_{n-1}x^{n-1}+\cdots+c_1x+c_0$ be a polynomial with rational coefficients.
%Suppose $f(x)$ cannot be factored as a product of polynomials with rational coefficients. If $f(\alpha)=0$ for some complex number $\alpha$, then we say that $\alpha$ is an algebraic number. Note that
%\[\alpha^n=-c_{n-1}\alpha^{n-1}-\cdots-c_1\alpha-c_0.\]
%Thus, the polynomial ring $\Q(\alpha)$ is isomorphic as an abelian group under addition to the vector space $\Q^n$.
%A field $K=\Q(\alpha)$, where $\alpha$ is an algebraic number is called a number field.
%A norm of an algebraic number $\beta$ is defined in the following way:
%if $g(x)$ is an irreducible polynomial over the rational numbers and $g(\beta)=0$ then the norm of $\beta$ is 
%\[N(\beta)=\prod_{i=1}^n\beta_i,\]
%where $\beta_1,\dots,\beta_n$ are all the roots of the polynomial $g(x)$.
%Let ${\mathcal{O}}_K\subset K$ be the ring of integers in $K$, which consists of all the elements $\beta$ of $K$ such that $g(\beta)=0$ for some (irreducible) polynomial with integer coefficients and leading term $1$. Let $U_K$ be the group of units in ${\mathcal{O}}_K$ which is the group of (multiplicatively) invertible elements in ${\mathcal{O}}_K$. Elements of $U$ have norm $\pm 1$.
%Then the Dedekind zeta function (for principal ideals) is defined as
%\[\zeta_K(s)=\sum_{\beta\in {\mathcal{O}}_K-\{0\}/U_K} \frac{1}{\left|N(\beta)\right|^s}.\]

%%%%%
\subsection{Multiple Dedekind zeta values (MDZVs)}

We recall the construction of MDZVs over a real quadratic field $K$. 
(See \cite{MDZF} for definition of MDZVs over any number field.)
Let ${\mathcal{O}}_K$ be the ring of integers in $K$.

And let  $\beta$ be a totally positive unit in  ${\mathcal{O}}_K$.
Let $C$ be the cone defined as $\N$-linear combination of $1$ and $\beta$, that is,
\[C=\{\gamma\in {\mathcal{O}}_K\,|\, \gamma=a+b\beta,\mbox{ for }a,b\in\N\}.\]
Let $f_0(C;t_1,t_2)=\sum_{\gamma\in C}\exp(-t_1\gamma_1-t_2\gamma_2),$
where $\gamma_1$ and $\gamma_2$ are two real embeddings of $\gamma$.
We express $\zeta_{K;C}(2)$, $\zeta_{K;C}(3)$ and $\zeta_{K;C}(1,2)$  as iterated integrals on a membrane.
See \cite{MDZF} and \cite{Hilbert}, for more examples and properties of iterated integrals on membranes.
\begin{eqnarray}
\label{int z2}
&&
\int_0^\infty\int_0^\infty\left(\int_{u_1}^\infty\int_{u_2}^\infty f_0(C;t_1,t_2)dt_1\wedge dt_2\right)du_1\wedge du_2\\
\nonumber
&&
=\int_0^\infty\int_0^\infty\left(\int_{u_1}^\infty\int_{u_2}^\infty  \left(\sum_{\gamma\in C}\exp(-t_1\gamma_1-t_2\gamma_2)\right)  dt_1\wedge dt_2\right)du_1\wedge du_2\\
\nonumber
&&
=\int_0^\infty\int_0^\infty \left(\sum_{\gamma\in C}\frac{\exp(-u_1\gamma_1-u_2\gamma_2)}{\gamma_1\gamma_2}\right) du_1\wedge du_2\\
\nonumber
&&=\sum_{\gamma\in C}\frac{1}{(\gamma_1\gamma_2)^2}
=\sum_{\gamma\in C}\frac{1}{N(\gamma)^2}=\\
\nonumber
&&=\zeta_{K;C}(2).
\end{eqnarray}

Similarly,
{\small \begin{multline*}
\zeta_{K;C}(3)
=\sum_{\gamma\in C}\frac{1}{N(\gamma)^3}\\
=\int_0^\infty\int_0^\infty
\left(\int_{v_1}^\infty\int_{v_2}^\infty
\left(\int_{u_1}^\infty\int_{u_2}^\infty 
f_0(C;t_1,t_2)dt_1\wedge dt_2\right)
du_1\wedge du_2\right)
dv_1\wedge dv_2,\\
\end{multline*}}
and
\begin{eqnarray*}
&&\zeta_{K;C}(1,2)
=\sum_{\gamma,\delta\in C}\frac{1}{N(\gamma)^1N(\gamma+\delta)^2}=\\
&&=\int_0^\infty\int_0^\infty
\bigg(\int_{v_1}^\infty\int_{v_2}^\infty
\bigg(\int_{u_1}^\infty\int_{u_2}^\infty
f_0(C;t_1,t_2)dt_1\wedge dt_2\bigg)\\
&&\times f_0(C;u_1,u_2)du_1\wedge du_2\bigg)
dv_1\wedge dv_2.
\end{eqnarray*}

%%%%%%%%%%%%%%%%%%%%%%%
\section{Transition to Algebraic Geometry}
We can write the infinite sum in the definition of $f_0$ as a product of two geometric series
\begin{multline*}
f_0(C;t_1,t_2)
=
\sum_{\gamma\in C}
\exp(-\gamma_1t_1-\gamma_2t_2)\\
=
\sum_{a=1}^\infty
\sum_{b=1}^\infty 
\exp[-(a\alpha_1+b\beta_1)t_1- (a\alpha_2+b\beta_2)t_2]\\
=
\sum_{a=1}^\infty
\sum_{b=1}^\infty 
\exp[-a(\alpha_1t_1+\alpha_2t_2)]
\exp[-b(\beta_1t_1+\beta_2t_2)]\\
=
\frac{\exp[-(\alpha_1t_1+\alpha_2t_2)]}{1-\exp[-(\alpha_1t_1+\alpha_2t_2)]}
\times
\frac{\exp[-(\beta_1t_1+\beta_2t_2)]}{1-\exp[-(\beta_1t_1+\beta_2t_2)]}
\end{multline*}
Let $x_1=e^{-t_1}$ and $x_2=e^{-t_2}$.
Then
\begin{equation}
\label{alg}
f_0(C;t_1,t_2)
=
\frac{x_1x_2}{1-x_1x_2}
\cdot
\frac{x_1^{\beta_1}x_2^{\beta_2}}{1-x_1^{\beta_1}x_2^{\beta_2}}
\end{equation}

Now we are going to express $f_0$ algebraically. At this point there is a problem of raising the variable $x$ to an integer algebraic power. Note that$\beta_1$ and $\beta_2$ are algebraic integers (in fact totally positive units), which are not rational integers.

How do we raise $x$ to power $\beta_1$ and to $\beta_2$? We introduce new variables 
\[y_1=x_1^{\beta_1}\mbox{ and }y_2=x_2^{\beta_2}.\]
Then $x_1^{a+b\beta_1}=x_1^{a}y_1^b$, where $a$ and $b$ are integers.

%\pagebreak

We are going to use the variables $x_1,x_2$. For each of then we introduce $y_1,y_2$, so that we write $y_1$ instead of 
$x_1^{\beta_1}$ and $y_2$ instead of $x_2^{\beta_2}$.
In terms of $x_1$, $x_2$, $y_1$ and $y_2$, we can express $f_0$ as
\begin{multline*}
f_0(C;t_1,t_2)=\frac{x_1x_2}{1-x_1x_2}\cdot\frac{x_1^{\beta_1}x_2^{\beta_2}}{1-x_1^{\beta_1}x_2^{\beta_2}}
=\frac{x_1x_2}{1-x_1x_2}\cdot\frac{y_1y_2}{1-y_1y_2}.
\end{multline*}

Let us also define  $\omega_1=\frac{d(x_1x_2)}{1-x_1x_2}\wedge\frac{d(y_1y_2)}{1-y_1y_2}$ and let
$\omega_0=\frac{d(x_1x_2)}{x_1x_2}\wedge\frac{d(y_1y_2)}{y_1y_2}$.

{\bf{Key Remark:}} The differential forms $\omega_0$ and $\omega_1$ will be used for both algebraic geometry on moduli spaces and for defining multiple Dedekind zeta values.

\begin{lemma}
If we substitute $x_1=e^{-t_1}$, $x_2=e^{-t_2}$, $y_1=e^{-\beta_1t_1}$ and $y_2=e^{-\beta_2t_2}$, then
\[\omega_0=(\beta_2-\beta_1)dt_1\wedge dt_2.\]
\end{lemma}
{\bf{Proof:}}
Consider $x_1,x_2,y_1$ and $y_2$ as functions of $t_1$ and $t_2$. Then
\[y_1y_2=x_1^{\beta_1}x_2^{\beta_2}\]
and
\[
\frac{d(y_1y_2)}{y_1y_2}=\frac{d(x_1^{\beta_1}x_2^{\beta_2})}{x_1^{\beta_1}x_2^{\beta_2}}=\beta_1\frac{dx_1}{x_1}+\beta_2\frac{dx_2}{x_2}=-\beta_1dt_1-\beta_2dt_2
\]
Similarly,
\[\frac{d(x_1x_2)}{x_1x_2}=-dt_1-dt_2.\]
Again, as functions of $t_1$ and $t_2$, we have
\begin{multline*}
\omega_0=\frac{d(x_1x_2)}{x_1x_2}\wedge\frac{d(y_1y_2)}{y_1y_2}
=
(dt_1+dt_2)\wedge(\beta_1dt_1+\beta_2dt_2)\\
=(\beta_2-\beta_1)dt_1\wedge dt_2.
\end{multline*}

Now let us write $\omega_0(x_1,x_2)$ and $\omega_1(x_1,x_2)$, when we want to specify the dependence on the variables.
In fact, both forms depend also on $y_1$ and $y_2$; however, we will take care of that by choosing a region of integration together with tangential base points.
%%%%%%%%%%%%%%%%%%
\section{Tangential base points}
Let $x_1=e^{-t_1}$ and let $y_1=e^{-\beta_1t_1}$ We would like to find  an algebraic relation among the variables $x_1$ and $y_1$ when they approach $(0,0)$ or when they approach $(1,1)$. That occurs when $t_1$ approaches $\infty$ or when $t_1$ approaches $0$, respectively. 
If $\beta_1>1$ then
\[\lim_{t_1\rightarrow \infty} \frac{dy_1}{dx_1}=
\lim_{t_1\rightarrow \infty} \frac{de^{-\beta_1t_1}}{de^{-t_1}}=
\lim_{t_1\rightarrow \infty}\beta_1 \frac{e^{t_1}}{(e^{ t_1})^{\beta_1}}=0.
\]
Also
\[\lim_{t_1\rightarrow 0} \frac{dy_1}{dx_1}=
\lim_{t_1\rightarrow 0} \beta_1\frac{e^{-\beta_1t_1}}{e^{-t_1}}=\beta_1.\]

Let \[\gamma_1:(0,\infty)\rightarrow \mathcal{M}_{0,5},\]
\[\gamma_1(t_1)=(e^{-t_1},e^{-\beta_1t_1})=(x_1,y_1).\]
For a vector $v=(a,b)$, consider $[v]=[a:b]$ as an element of $\P^1$.

We have proven the following lemma.
\begin{lemma}

(a)  $\lim_{t_1 \rightarrow \infty}\left[\frac{d\gamma_1}{dt_1}\right]=[1:0],$

(b) $\lim_{t_1 \rightarrow 0}\left[\frac{d\gamma_1}{dt_1}\right]=[1:\beta_1].$
\end{lemma}

Similarly,
we have $x_2=e^{-t_2}$ and $y_2=e^{-\beta_2t_2}$ with $0<\beta_2<1$.
Let \[\gamma_2:(0,\infty)\rightarrow \mathcal{M}_{0,5},\]
\[\gamma_2(t_2)=(e^{-t_2},e^{-\beta_2t_2})=(x_2,y_2).\]
The following Lemma could be proven in the same way.

\begin{lemma}

(a) $\lim_{t_2 \rightarrow \infty}\left[\frac{d\gamma_2}{dt_2}\right]=[0:1],$

(b) $\lim_{t_2 \rightarrow 0}\left[\frac{d\gamma_1}{dt_2}\right]=[1:\beta_2].$
\end{lemma}

{\bf{Remark:}} The paths $\gamma_1$ and $\gamma_2$ can be used to define a membrane $m=\gamma_1\times\gamma_2$ by taking a Cartesian products of both the domains and the targets 
\[
m=\gamma_1\times\gamma_2:(0,1)^2 \rightarrow \left(\mathcal{M}_{0,5}\right)^2.
\] 
The definition of multiple Dedekind zeta values via iterated integrals on a membrane use exactly the membrane $m$ in the case of quadratic fields (see \cite{MDZF}).

\begin{proposition} With the above choice of tangential base points, we have
\[\int_{0<x_1<x_3<1;\,\,0<x_2<x_4<1}\omega_1(x_1,x_2)\wedge \omega_0(x_3,x_4)=(\beta_2-\beta_1)^2\zeta_{K;C}(2).\]
\end{proposition}
{\bf{Proof:}}
The differential forms $\omega_0$ and $\omega_1$  are closed. Thus we can vary the paths $\gamma_1$ and $\gamma_2$ without changing the value of the integral as long as the tangential base points remain the same. Thus, we can choose the parametrization
$x_i=e^{-t_i}$ and $y_i=e^{-\beta_it_i}$, keeping the tangential points fixed. Using Formulas \eqref{int z2} and \eqref{alg}, 
we obtain

\begin{eqnarray*}
\frac{d(x_3x_4)}{x_3x_4}
\wedge
\frac{d(y_3y_4)}{y_3y_4}
=
(\beta_2-\beta_1)dt_3\wedge dt_4
\end{eqnarray*}
Similarly,
we have that

\begin{multline*}
\frac{x_1x_2}{1-x_1x_2}
\cdot
\frac{y_1y_2}{1-y_1y_2}
\cdot
\left(
\frac{d(x_3x_4)}{x_3x_4}
\wedge
\frac{d(y_3y_4)}{y_3y_4}
\right)\\
=
f_0(C;t_1,t_2)(\beta_2-\beta_1)dt_1\wedge dt_2.
\end{multline*}
Thus, with the above choice of tangential base points, we have
\begin{eqnarray*}
&&\int_{0<x_1<x_3<1;\,\,0<x_2<x_4<1}\omega_1(x_1,x_2)\wedge \omega_0(x_3,x_4)\\
&&=(\beta_2-\beta_1)^2\int_{t_1>t_3>0;\,\,t_2>t_4>0}f_0(C;t_1,t_2)dt_1\wedge dt_2 \wedge t_3 \wedge dt_4\\
&&=(\beta_2-\beta_1)^2\zeta_{K;C}(2).
\end{eqnarray*}

\begin{corollary} With the above choice of tangential base points, we have
\label{cor}

%(a) 
%\begin{multline*}
%(\beta_2-\beta_1)^n\zeta_{K;C}(n)\\
%=\int_{ 0<x_1<x_3<\dots<x_{2n-1}<1;\,\,0<x_2<x_4<\dots<x_{2n}<1}
%\omega_1(x_1,x_2)\\
%\wedge \omega_0(x_3,x_4)\wedge\dots\wedge \omega_0(x_{2n-1},x_{2n})
%\end{multline*}
%(b)
\begin{multline*}
(\beta_2-\beta_1)^3\zeta_{K;C}(1,2)\\
=
\int_{0<x_1<x_3<x_{5}<1;\,\,0<x_2<x_4<x_{6}<1}
\omega_1(x_1,x_2)\wedge \omega_1(x_3,x_4)\wedge \omega_0(x_{5},x_{6}).
\end{multline*}
\end{corollary}

\begin{theorem} In Corollary \ref{cor}, the integral on the right hand side is  a period of a mixed Tate motive unramified over a real quadratic number ring.
\end{theorem}
{\bf{Proof:}} In this proof we are going to follow closely the paper by Goncharov and Manin \cite{GM}.
The period will be a pairing between 
 $[\Omega_A]\in Gr^W_{12}H^6(\overline{\mathcal{M}}_{0,15}-A)$ and $[\Delta_B]\in \left(Gr^W_0H^6(\overline{\mathcal{M}}_{0,15}-B)\right)^\vee$ associated to a mixed Tate motive $H^6(\overline{\mathcal{M}}_{0,15}-A;B-A\cap B)$.
 
Let the $(4n)$-coordinates $x_{2i-1},y_{2i-1},z_{2i-1},w_{2i-1}$ for indices $i=1,2,\dots,n,$ be a coordinate of a point on $\mathcal{M}_{0,4n+3}$.
One can think of  $\mathcal{M}_{0,4n+3}$ as $(\P^1)^{4n}-D$ where the divisor $D$ is obtained by setting any of the coordinates to be $0$, $1$, $\infty$  or setting any two of the coordinates to be equal. Let us define 
\[x_{2i}=\frac{1}{z_{2i-1}}\mbox{ and }y_{2i}=\frac{1}{w_{2i-1}}.\]
Now the coordinates of any point on $\mathcal{M}_{0,4n+3}$ can be written as $(x_1,y_1,x_2,y_2,\dots,x_{2n},y_{2n})$. In terms of the new coordinates, we have the following components of $D$:

$x_i=0$, $x_i=1$, $x_i=\infty$,

$y_i=0$, $y_i=1$, $y_i=\infty$,

$x_1=x_3$, $x_3=x_5$,

$y_1=y_3$, $y_3=y_5$,

$x_1x_2=1$, 

$x_3x_4=1$,

$y_1y_2=1$, 

$y_3y_4=1$.

The last four components can be realized in terms of the previous coordinates as $x_1=z_1$, $x_3=z_3$, $y_1=w_1$ and $y_3=w_3$.

Let $n=3$. Let $\overline{\mathcal{M}}_{0,4n+3}=\overline{\mathcal{M}}_{0,15}$ be the Deligne-Mumford compactification of the moduli space of curves of genus $0$ with $15$ marked points. The ambient space will be $\overline{\mathcal{M}}_{0,15}$. From it we will remove a divisor $A$ whose components occur as poles of the differential forms under the integral. Explicitly, the differential forms are
\[\omega_1(x_1,x_2)=\frac{d(x_1x_2)}{1-x_1x_2}\wedge\frac{d(y_1y_2)}{1-y_1y_2},\]
\[\omega_1(x_3,x_4)=\frac{d(x_3x_4)}{1-x_3x_4}\wedge\frac{d(y_3y_4)}{1-y_3y_4},\]
\[\omega_0(x_{5},x_{6})=\frac{d(x_5x_6)}{x_5x_6}\wedge\frac{d(y_5y_6)}{y_5y_6}.\]

The components of the divisor $A$ consists of the uinion of

$(x_1x_2=1)$, $(y_1y_2=1)$, $(x_3x_4=1)$, $(y_3y_4=1)$,

$(x_5=0)$,  $(x_6=0)$,  $(y_5=0)$,  $(y_6=0)$, 
 
$(x_i=\infty)$,   $(y_i=\infty)$, for $i=1,2,\dots,6$,
  
together with the exceptional divisors obtained via blow-up at the intersections of two components that both contain the same variable or the same constant $0$, $1$ or $\infty$ on the right hand side of the equalities.

Thus, the differential form
\[\Omega_A=\omega_1(x_1,x_2)\wedge \omega_1(x_3,x_4)\wedge \omega_0(x_{5},x_{6})\]
is well-defined on $\overline{\mathcal{M}}_{0,15}-A$.

Now we proceed to defining $B$. The key part will be to include the tangential base points in the definition of $B$.

The components of $B$ consist of a union of codimension $1$ subvarieties and codimension $2$ subvarieties. The latter ones correspond to the tangential base points.  

The codimension $1$ components are the following:

$(x_1=0)$, $(x_1=x_3)$, $(x_3=x_5)$, $(x_5=1)$,

$(x_2=0)$, $(x_2=x_4)$, $(x_4=x_6)$, $(x_6=1)$,

$(y_1=0)$, $(y_1=y_3)$, $(y_3=y_5)$, $(y_5=1)$,

$(y_2=0)$, $(y_2=y_4)$, $(y_4=y_6)$, $(y_6=1)$,

together with the exceptional divisors of the blow-up at an intersection of two subvarieties such that the two polynomials contain the same variable or the same constant $0$ or $1$ on the right hand side of the equaities, except the following $4$ double intersections of components

$(x_1=0)$ and $(y_1=0)$,

$(x_2=0)$ and $(y_2=0)$,

$(x_5=1)$ and $(y_5=1)$,
 
$(x_6=1)$ and  $(y_6=1)$,

to which we associate a codimension $2$ subvarieties of $\overline{\mathcal{M}}_{0,15}$, using the tangential base points.

For the blow-up at the intersection $(x_1=0)$ and $(y_1=0)$ we choose a divisor $B_1$ on the exceptional divisor defined by $[x_1:y_1]=[1:0]$. 
Note that $B_1$ is of codimension $2$ in $\overline{\mathcal{M}}_{0,15}$.

For the blow-up at the intersection $(x_2=0)$ and $(=0y_2)$ we choose a divisor $B_2$ on the exceptional divisor defined by $[x_2:y_2]=[0:1]$. 

For the blow-up at the intersection $(x_5=1)$ and $(y_5=1)$ we choose a divisor $B_5$ on the exceptional divisor defined by $[x_5:y_5]=[1:\beta_1]$. 

For the blow-up at the intersection $(x_6=1)$ and $(y_6=1)$ we choose a divisor $B_6$ on the exceptional divisor defined by $[x_6:y_6]=[1:\beta_2]$. 

The tangential base points define the components $B_1,B_2,B_5,B_6$. Thus,
$(\beta_2-\beta_1)^3\zeta_{K,C}(1,2)$ occurs as a period of
$H^6(\overline{\mathcal{M}}_{0,15}-A;B-A\cap B)$ when $[\Omega_A]\in Gr^W_{12}H^{6}(\overline{\mathcal{M}}_{0,15}-A)$ is paired with 
$[\Delta_B]\in \left(Gr^W_0H^6(\overline{\mathcal{M}}_{0,15}-B)\right)^{\vee}$

Note that $B_1$ and $B_2$ are defined over $\Z$, and $B_5$ and $B_6$ are defined over the ring of integers $\mathcal{O}_K$ of the field $K$. Each of them is naturally isomorphic to $\overline{\mathcal{M}}_{0,13}$ as a variety over $\mathcal{O}_K$. 
Similarly, any intersection of the components of $B$ is isomorphic  over $\mathcal{O}_K$ to $\overline{\mathcal{M}}_{0,n}$ for some integer $n$.
Using that $H^i(\overline{\mathcal{M}}_{0,n})$ is a mixed Tate motive over $Spec(\mathcal{O}_K)$, we obtain that the motivic cohomology of the components of $B$ are mixed Tate motives. Using Proposition 1.7 from Deligne and Goncharov, \cite{DG}, we conclude that for $l\neq char (\nu)$ the $l$-adic cohomology of the reduction of $B_j$  modulo $\nu$ of the motive $H^i(B_j)$ is unramified for any component $B_j$ of $B$, since $B_j$ is isomorphic to $\overline{\mathcal{M}}_{0,n}$ over $Spec(\mathcal{O}_K)$  for some $n$. We conclude that for $l\neq char (\nu)$ the $l$-adic cohomology of the reduction modulo any $\nu\in Spec(\mathcal{O}_K)$ of the motive $H^6(\overline{\mathcal{M}}_{0,15}-A;B-A\cap B)$ is unramified. Thus,
$H^6(\overline{\mathcal{M}}_{0,15}-A;B-A\cap B)$ is a mixed Tate motive unramified over $Spec(\mathcal{O}_K)$.

%%%%%%%%%%%%%%%%%%%%%%%%%%%%%%%%%%%%%%%%%%%%
\section*{Acknolwedgements}
I would like to thank the organizers of the XXXV Workshop on
Geometry and Mathematical Physics in Biaowieza, Poland, especially, Emma Previato, for the opportunity to present this result, and for the encouragement and interest in my work. I would also like to thank the referees of the paper for their questions, suggestions and corrections. Due to their comments, the mathematical content of the paper has been improved immensely. 

%Due to their work, the paper is substantially improved - both stylistically and mathematically. I would like to express special thanks to the referees for their questions and suggestions regarding the proof of the Theorem.

%\renewcommand{\em}{\textrm}

%\begin{small}

%\renewcommand{\refname}{ {\flushleft\normalsize\bf{References}} }

%\end{small}

%Ivan Horozov

City University of New York,
Bronx Community College
2155 University Avenue, Bronx,
New York 10453 , U.S.A.;
ivan.horozov@bcc.cuny.edu

\end{document}